\documentclass[11pt]{amsart}
\usepackage{epsfig,amsmath,amssymb,amsthm,epsf}
\usepackage{graphicx,array}

\newcommand{\zz}[1]{\mathbb #1}
\newtheorem{proposition}{Proposition}
\newtheorem{lemma}{Lemma}
\newtheorem{conjecture}{Conjecture}
\newtheorem{theorem}{Theorem}
\newtheorem{corollary}{Corollary}

\begin{document}
\title{{Ergodicity and Mixing Properties of the Northeast Model }}
\author{{George Kordzakhia}}
\address{University of California\\
Department of Statistics\\
Berkeley CA }
\email{kordzakh@stat.berkeley.edu}
\author{Steven P. Lalley} 
\address{University of Chicago\\ Department of Statistics \\ 5734
University Avenue \\
Chicago IL 60637}
\email{lalley@galton.uchicago.edu}
\date{\today}
\maketitle
\begin{abstract}
The \emph{northeast model} is a spin system on the two-dimensional
integer lattice that evolves according to the following rule: Whenever
a site's southerly and westerly nearest neighbors have spin $1$, it
may reset its own spin by tossing a $p-$coin; at all other times, its
spin remains frozen. It is proved that the northeast model has a phase
transition at $p_{c}=1-\beta_{c}$, where $\beta_{c}$ is the critical
parameter for oriented percolation.  For $p<p_{c}$, the trivial
measure $\delta_{0}$ that puts mass one on the configuration with all
spins set at $0$ is the unique ergodic, translation invariant,
stationary measure. For $p\geq p_{c}$, the product Bernoulli-$p$
measure on configuration space is the unique nontrivial, ergodic,
translation invariant, stationary measure for the system, and it is
mixing. For $p>2/3$ it is shown that there is exponential decay of
correlations. 
\end{abstract}


\section{Introduction}\label{sec:intro} The \emph{northeast model} is
the simplest nontrivial facilitated spin-flip system on the
two-dimensional integer lattice $\zz{Z}^{2}$. Interest in such systems
dates to \cite{FA84}; for a recent review see \cite{PYA00}. The
northeast model is specified by the following rules: Spins take values
in the two-element set $\{0,1\}$; the spin at site $(x,y)\in
\zz{Z}^{2}$ may flip only at times when the spins at sites $(x-1,y)$
and $(x,y-1)$ are both $1$; and at such times the flip rates are
\begin{align}\label{eq:NE}
0 & \longrightarrow 1 & \quad &\text{at rate} \; p;\\
\notag 
1 & \longrightarrow 0 & \quad &\text{at rate} \; 1-p.
\end{align}
A rigorous construction of the process 
%using a system of marked Poisson processes 
is outlined in sec.~\ref{sec:bg}. It will also be shown that the
product Bernoulli-$p$ measure $\nu_{p}$ on the configuration space is
invariant and reversible. Note that unlike many growth models (e.g.,
the contact process) the northeast model is \emph{not} additive (in
the sense of Liggett \cite{liggett}, Definition III.6.1) or even
monotone (in the sense that stochastic monotonicity is preserved by
the evolution).

The northeast model is the natural two-dimensional analogue of the
one-dimensional \emph{east model} studied by Aldous and Diaconis
\cite{aldous-diaconis}. The east model is the spin-flip model on the
integer lattice $\zz{Z}$ in which the spin at site $x$ may flip only
at times when the spin at site $x-1$ is $1$; the flip rates at such
times are given by (\ref{eq:NE}).  The main results of
\cite{aldous-diaconis} are rigorous bounds on relaxation rates to the
product Bernoulli-$p$ measure that are sharp in the limit $p
\rightarrow 0$. The northeast model differs fundamentally from the
east model in that the product Bernoulli-$p$ measure $\nu_{p}$ is not
even ergodic for $p<1-\beta_{c}$, where $\beta_{c}$ is the critical
parameter for oriented site percolation. The reason for this is
apparent: If $p<1-\beta_{c}$ then $\nu_{p}-$almost every spin
configuration must contain infinite southwest clusters of 0s; these
clusters must remain frozen at spin $0$ forever.

We shall establish in this paper that the northeast model
exhibits a phase transition at $p_{c}:=1-\beta_{c}$  by showing that
for $p>p_{c}$ the product Bernoulli-$p$ measure is ergodic and mixing:

\begin{theorem}\label{theorem:1}
The only ergodic, translation invariant, stationary measures for the
northeast model with flip parameter $p\geq 1-\beta_{c}$ are $\nu_{p}$
and $\nu_{0}$. The only ergodic, translation invariant, stationary
measure for the northeast model with flip parameter $p<1-\beta_{c}$ is
$\nu_{0}$. For $p>p_{c}$, the product Bernoulli measure $\nu_{p}$  is
not only ergodic but also mixing.
\end{theorem}

Here \emph{mixing} means that correlations decay to zero as $t
\rightarrow \infty$. Following is a more precise statement of what we
shall prove. For any spin configuration $\zeta $, denote by $\zeta^{\Lambda}$ its
restriction to the set $\Lambda \subset \zz{Z}^{2}$. Similarly, denote
by $\nu_{p}^{\Lambda}$ the product Bernoulli-$p$ measure on the
restricted configuration space $\{0,1 \}^{\Lambda}$.  Denote by
$\xi_{t}$ the (random) spin configuration at time $t$ in a realization
of the northeast model.

\begin{theorem}\label{theorem:2}
If $p>1-\beta_{c}$ then $\nu_{p}$ is mixing for the northeast model in
the following sense: For $\nu_{p}-$almost every initial configuration
$\sigma$, and for every finite set $\Lambda \subset \zz{Z}^{2}$ and
every configuration $\zeta^{\Lambda}\in \{0,1 \}^{\Lambda}$, 
\begin{equation}\label{eq:mixing}
	\lim_{t \rightarrow \infty}
	 P_{p}^{\sigma} \{\xi_{t}^{\Lambda}=\zeta^{\Lambda } \}
	=\nu_{p}^{\Lambda} (\zeta^{\Lambda}).
\end{equation}
\end{theorem}

In light of the Aldous-Diaconis results, it is
natural to ask if correlations decay exponentially for $p>p_{c}$. We
conjecture that they do. We have been able to prove this only for
$p>2/3$. 

\begin{theorem}\label{theorem:3}
For each $p>2/3$ and each finite set $\Lambda$ of sites there exist
constants $C_{\Lambda},\alpha_{\Lambda}>0$ such that  for any two
configurations $\zeta^{\Lambda}, \eta^{\Lambda}$, $\forall \, t>0$,
\begin{equation}\label{eq:expMixing}
	|P_{p}^{\nu_{p}} \{\xi_{t}^{\Lambda}=\zeta^{\Lambda };
\xi_{0}=\eta^{\Lambda} \}- \nu_{p}^{\Lambda}(\zeta^{\Lambda})\nu_{p}
(\eta^{\Lambda}) | \leq C_{\Lambda}\exp \{-\alpha_{\Lambda}t \}.
\end{equation}
\end{theorem}

We conjecture further that the exponential decay parameters
$\alpha_{\Lambda}=\alpha_{\Lambda} (p)$ can be chosen so as not to
depend on the sets $\Lambda$. If this is true then for $p>p_{c}$ the
northeast model has a positive \emph{spectral gap} (see
\cite{aldous-diaconis}), and the nature of the phase transition would
be reflected in the behavior of the spectral gap as $p \rightarrow
p_{c}$.

Sections \ref{sec:bg}--\ref{sec:expMixing} are devoted to the proofs
of Theorems \ref{theorem:1}--\ref{theorem:3}. Theorems \ref{theorem:1}
and \ref{theorem:2} will be proved in
secs.~\ref{sec:bg}--\ref{sec:ergodicity}, and Theorem \ref{theorem:3}
in sec.~\ref{sec:expMixing}. In the final section \ref{sec:shape} we
state a conjecture about the propagation of influence in the northeast
model.

\section{Construction}\label{sec:bg} Two constructions of the
northeast model are possible, one using the Hille-Yosida theorem (see
Liggett \cite{liggett}, ch.~1), the second using a system of
independent marked Poisson processes. The first has the advantage that
it gives simple characterizations of invariant and reversible
measures. The second yields detailed information about the time
evolution of the process.

\subsection{Generator and Semigroup}\label{ssec:generator}
Let $X$ be the space of spin configurations on the lattice
$\zz{Z}^{2}$, and for any subset $\Lambda \subset \zz{Z}^{2}$ let
$X^{\Lambda}$ be the set of configurations on
$\Lambda$. Denote by $C (X)$ the space of continuous real-valued
functions  on $X$ (relative to the product topology), and by
$C_{*}(X)$ the subset consisting of those functions that depend only
on finitely many coordinates. Note that $C_{*} (X)$ may be naturally
identified with $\cup_{\Lambda} C (X^{\Lambda})$, where the  union is
over all finite subsets $\Lambda $ of $\zz{Z}_{2}$. For any
configuration $\sigma$, define $A (\sigma )$ to be the set of sites
that are flip-eligible in configuration $\sigma$, that is, those sites
$x$ whose nearest neighbors to the south and west both have spin $1$
in configuration $\sigma$. Also, for any site $x$, denote by $\sigma^{+}_{x}$
and $\sigma^{-}_{x}$ the configurations that agree with $\sigma$ at
all sites $y\not =x$ and have values $1,0$, respectively, at site
$x$. 

For any function
$f\in C (X)$ and each site $x\in \zz{Z}_{2}$ define
\begin{align}\label{definition:del}
	\nabla_{x}f (\sigma )&= ( f (\sigma^{+}_{x}) -f
	 (\sigma^{-}_{x}))p &\quad &\text{if} \; \sigma =\sigma^{-}_{x};\\
\notag 	&=( f (\sigma^{-}_{x}) -f
	 (\sigma^{+}_{x})) (1-p)  &\quad &\text{if} \; \sigma =\sigma^{+}_{x}.
\end{align}
For any $f\in C_{*} (X)$ define
\begin{equation}\label{definition:generator}
	\mathcal{L}f (\sigma ) 
	 =\sum_{x\in \zz{Z}_{2}} \nabla_{x}f (\sigma)
	  \mathbf{1}\{x\in A (\sigma)\};
\end{equation}
Observe that because $f\in C_{*}
(X)$ depends on only finitely many coordinates the sum is finite, and
hence the definition is a valid one. It follows from \cite{liggett},
sec.~1.3 that the operator $\mathcal{L}$ on $C_{*} (X)$ extends
uniquely to a Markov generator, and that $C_{*} (X)$ is a core for the
generator. (Note: \cite{liggett} uses a larger core, but the graph of
$\mathcal{L}$ restricted to $C_{*} (X)$ is easily seen to have closure
equal to that of $\mathcal{L}$ restricted to the larger core of
\cite{liggett}.) The Hille-Yosida theorem therefore implies that the
closure of $\mathcal{L}$ generates a unique Feller Markov semigroup on
$C (X)$. This implies the existence of a Markov process with flip
rates (\ref{eq:NE}).

\subsection{Stationary and Reversible Measures}\label{ssec:invariance}
A Borel probability measure $\mu$ on $X$ is stationary for the Markov
semigroup generated by $\mathcal{L}$ if and only if for every $f\in
C_{*} (X)$,
\begin{equation}\label{eq:stationarity}
	\int \mathcal{L}f \,d\mu =0.
\end{equation}
A stationary measure $\mu$ is reversible if and only if for all pairs
$f,g \in C_{*} (X)$,
\begin{equation}\label{eq:reversibility}
	\int f\mathcal{L}g \,d\mu =
	\int g\mathcal{L}f \,d \mu .
\end{equation}
These criteria yield easy proofs that the product Bernoulli-$p$
measure $\nu_{p}$ is stationary and reversible for the northeast model
with flip parameter $p$. By (\ref{definition:generator}) and
(\ref{eq:stationarity}), to prove that $\nu_{p}$ is stationary it
suffices to show that for any function $f:X \rightarrow \zz{R}$ that
depends only on the spins in $\Lambda$ and each $x\in \Lambda$,
\[
	\int_{\{\sigma \in F (x) \}}
	 \nabla_{x}f (\sigma)\,
	 d\nu_{p} (\sigma)=0
\]
where $F (x)=\{\sigma : x\in A (\sigma) \}$ is the set of
configurations $\sigma$ for which site $x$ is flip-eligible.
The event $\{ x\in A (\sigma)\}=\{\sigma \in F (x) \}$ does not
depend on the spin at $x$, and so under $\nu_{p}$ this event is
independent of $\sigma (x)$. Moreover, the event $\sigma
=\sigma^{+}_{x}$ occurs if and only if $\sigma (x)=1$, and hence is
independent of $\{\sigma \in F (x) \}$. Since $\sigma (x)=1$ with
probability $p$ under $\nu_{p}$, (\ref{definition:del}) implies that
(with $q=1-p$)
\begin{equation*}
\begin{split}
	&\int_{F (x)} \nabla_{x}f (\sigma) \, 
	 d\nu_{p} (\sigma)\\
	&\quad =\int_{F (x)} ((f (\sigma^{+}_{x})-f (\sigma^{-}_{x}))pq 
	  +(f (\sigma^{-}_{x})-f (\sigma^{+}_{x}))qp )\, d\nu_{p} (\sigma)\\
	&\quad =0.
\end{split} 
\end{equation*}
Similarly, to prove that $\nu_{p}$ is reversible, it suffices to show
that for any pair of functions $f,g \in C (X)$ that only depend on the
spins in a finite set $\Lambda$,
\[
	\int_{F (x)} g (\sigma )\nabla_{x}f (\sigma) \,
	 d\nu_{p} (\sigma)=
	\int_{F (x)} f(\sigma )\nabla_{x}g (\sigma) \,
	 d\nu_{p} (\sigma).
\]
Since the spin at $x$ is independent of the event $\sigma \in F (x)$,
this identity follows by a simple calculation similar to that above.

\subsection{Construction via Marked Poisson Processes}\label{ssec:ppConstruction}

Let $(\Omega ,\mathcal{F},P)$ be a probability space on which are
defined countably many independent rate-one Poisson processes and
countably many independent uniform-(0,1) random variables. Assign one
Poisson process to each site $x\in \zz{Z}^{2}$, and one uniform
r.v. to each occurrence in each Poisson process. The occurrences in
the Poisson process attached to site $x$ mark the times of flip
opportunities for site $x$: at each such time, site $x$ queries its
neighbors to the south and west about their current spins, and if both
of these are $1$ then $x$ resets its spin according to the value $U$
of the uniform r.v. attached to the occurrence. The reset rule is
\begin{align}\label{eq:spinDetermination}
	U&\leq p & \Longrightarrow \quad  &\text{reset spin to } \;1;\\
\notag 	U&> p    & \Longrightarrow \quad &\text{reset spin to } \;0.
\end{align}
Note that such a reset does not necessarily change the spin at $x$: if
the reset occurs at a time when the spin at $x$ is $0$ (respectively,
$1$), then the chance that the spin changes is $p$ (respectively,
$q$). Thus, the flip rates at site $x$ will agree with the
specification (\ref{eq:NE}).

In order that this construction uniquely specify the spin at site $x$
at each finite time, it must be the case that the spins of the
southwest neighbors of $x$ at query times can be determined. For this
to be the case, the initial configuration  $\xi_{0}$ must be
specified. Throughout the remainder of this section and section
\ref{ssec:regeneration} we assume that $\xi_{0}=\zeta$ for some fixed,
nonrandom configuration $\zeta$. (We suppress the dependence of the
process $\xi_{t}$ on $\zeta$, but the reader should note that
the algorithm specified below will give a different time evolution
$\xi_{t}$ for every different choice of $\zeta$.)  To
determine the spin at site $y$ at time $t$, one must search backward
in time for occurrences in the Poisson process attached to site $y$:
If there are no occurrences in this process between times $0$ and $t$,
then the spin at site $y$ at time $t$ is set to its value at time
$0$. Otherwise, at the last occurrence time $s$ before $t$, the
backward search moves to the south and west neighbors $y',y''$ of $y$
and proceeds recursively. If it can then be determined that the spins
at $y'$ and $y''$ at time $s$ are both $1$ then the spin at $y$ at
time $t$ is determined by the attached uniform according to the rule
(\ref{eq:spinDetermination}). Otherwise, the query continues to the
last occurrence in the Poisson process at $y$ before time $s$. To show
that this algorithm terminates, we must show that the backward tree of
queries initiated at $(y,t)$ is almost surely finite. But this follows
because the backward query tree is stochastically dominated by a
simple binary fission process with fission rate $1$. (At each time a
query is made, two new query processes are engendered. The additional
queries in these offspring processes are mutualy independent,
\emph{except} when they coincide by virtue of a merger, e.g., SW and
WS).

\subsection{Regeneration at Reset Times}\label{ssec:regeneration}

Assume henceforth that the northeast process $\xi_{t}$ has been
constructed as in sec.~\ref{ssec:ppConstruction}, using auxiliary
marked Poisson processes and a specified initial configuration
$\xi_{0}$. For each site $x$, define 
\begin{equation}\label{definition:tau}
	\tau_{x}= \;\text{time of first reset at } \; x.
\end{equation}
Note that $\tau_{x}$ may take the value $+\infty$. Note also that
$\tau_{x}$ need not be (and generally will not be) measurable with
respect to the $\sigma -$algebra generated by $\{\xi_{t} \}_{0\leq
t<\infty}$.

Let $\leq$ and $<$ be the natural weak and strong partial orders on
$\zz{Z}^{2}$: that is, $x\leq y$ (respectively, $x<y$) if each
coordinate of $x$ is $\leq$ (respectively, $<$)
the corresponding coordinate of $y$. For each site $x$, define
\begin{align}\label{definition:ideals}
\mathcal{A}_{x}&=\{y \, : \, y \leq x \} 
\quad \text{and}\\
\notag 
\mathcal{B}_{x}&=\{y \, : \, x\not \leq y\}. &
\end{align}
Sites in $\mathcal{A}_{x}$ are those that may influence the times at
which spin resets at $x$ may occur; and sites in $\mathcal{B}_{x}$ are
those whose reset times are not influenced by site $x$. Let
$N^{x}_{t}$ be the the number of occurrences in the Poisson process
attached to site $x$ up to time $t$, and let $\{ U^{x}_{j}\}_{j\geq
1}$ be the uniform-$(0,1)$ random variables attached to the
occurrences of this Poisson process, listed in chronological
order. Define $\sigma -$algebras
\begin{align}\label{definition:sigmas}
\mathcal{G}^{x}_{t}&=\sigma (N^{x}_{s})_{s\leq t};\\
\notag 
\mathcal{H}^{x}_{t}&=\sigma (U^{x}_{i})_{i\leq N^{x}_{t}}; \\
\notag 
\mathcal{F}^{x}_{t}&=\sigma (\mathcal{G}^{x}_{t}\cup \mathcal{H}^{x}_{t});\\
\notag 
\mathcal{I}^{x}_{t}&= \sigma (\cup_{y\in \mathcal{A}_{x}}
	\mathcal{F}^{y}_{t});  
\quad \\
\notag 
\mathcal{I}^{\infty }_{t}&= \sigma (\cup_{x}\mathcal{I})^{x}_{t});
\\
\notag 
\mathcal{J}^{x}_{t}&=\sigma (\cup_{y\in
\mathcal{B}_{x}}\mathcal{F}^{y}_{t}); \quad \text{and}\\
\notag 
\mathcal{K}^{x}_{t}&=\sigma (\mathcal{J}^{x}_{t}\cup \mathcal{G}^{x}_{t}).
\end{align}
Note that $\mathcal{I}^{\infty}_{t}$ is generated by everything that
happens in the entire construction up to time $t$.

\begin{proposition}\label{proposition:resets}
For each site $x$ and each time $t>0$, the conditional distribution of
the spin $\xi_{t} (x)$ given  the
$\sigma -$algebra $\mathcal{K}^{x}_{\infty}$ on the event
$\{\tau_{x}\leq t \}$ is Bernoulli-$p$, that is,
\begin{equation}\label{eq:resets}
	P (\xi_{t} (x)=1 \, | \,\mathcal{K}^{x}_{\infty})=p
		\qquad \text{on} \;\{\tau_{x}\leq t \}.
\end{equation}
\end{proposition}

\begin{proof}
Observe that the reset time $\tau_{x}$, and in fact the entire string
of flip opportunity times at site $x$ up to time $t$, are measurable
with respect to the $\sigma -$algebra $\mathcal{K}^{x}_{t}$. On the
event $\tau_{x}\leq t$, there is a last flip opportunity at $x$ before
time $t$, and this occurs at the $k$th occurrence of $N^{x}_{s}$ for
some $k\geq 1$. The spin at $x$ at time $t$ is determined by the
uniform r.v. $U^{x}_{k}$. But this uniform is independent of
$\mathcal{K}^{x}_{\infty}$. The assertion (\ref{eq:resets}) follows.
\end{proof}

\begin{corollary}\label{corollary:mixing}
If,  for every site $x$ in a finite set $\Lambda$,
\begin{equation}\label{eq:tauFinite}
	P\{\tau_{x}<\infty \}=1
\end{equation}
then the joint distribution of the restricted configuration
$\xi^{\Lambda}_{t}$ converges weakly to the product Bernoulli-$p$
measure $\nu^{\Lambda}_{p}$ as $t \rightarrow \infty$.
\end{corollary}

\begin{proof}
Consider first for illustration the case $\Lambda =\{x,y \}$ where $x$
is the westerly nearest neighbor of $y$. Observe that
$\mathcal{K}^{x}_{t}\subset \mathcal{K}^{y}_{t}$, and recall that
$\xi_{t} (x)$ is measurable with respect to $\mathcal{K}^{y}_{t}$. By
hypothesis, $P\{\tau_{x}\leq t \}\approx 1$ and $P\{\tau_{y}\leq t
\}\approx 1$ for large $t$. Consequently, by Proposition
\ref{proposition:resets} (first for $y$, then $x$), for large $t$,
\begin{equation*}
\begin{split} 
	P (\xi_{t} (x)=\xi_{t} (y)=1) 
&	\quad =E E(P (\xi_{t} (y)=1 \, | \, \mathcal{K}^{y}_{t})
	 	\xi_{t} (x) \, | \,\mathcal{K}^{x}_{t})\\
&	\quad \approx p E P (\xi_{t} (x)=1 \, | \, \mathcal{K}^{x}_{t})\\
&	\quad \approx p^{2}.
\end{split} 
\end{equation*}

The general case is proved by induction on the size of
$\Lambda$. Fix a site $y\in \Lambda $ so that $\Lambda \setminus \{y
\}\subset \mathcal{B}_{y}$: since $\Lambda $ is finite, it must
contain such a site.  Set $\Lambda '=\Lambda \setminus
\{y\}$. Observe that for each $x\in \Lambda '$, the spin
$\xi_{t} (x)$ is measurable with respect to $\mathcal{K}^{y}_{t}$. Let
$\zeta^{\Lambda}$ be any spin configuration on $\Lambda$, and let
$\zeta^{R}$ be its restriction to a subset $R$ of $\Lambda$. Then by
Proposition \ref{proposition:resets} and the induction hypothesis, for
large $t$,
\begin{equation*}
\begin{split} 
	P (\xi^{\Lambda}_{t}=\zeta^{\Lambda})
&	=E P (\xi_{t} (y)=\zeta^{y} \, | \, \mathcal{K}^{y}_{t})
 	\mathbf{1}\{\xi^{\Lambda '}_{t}=\zeta^{\Lambda '} \} \\
&	\approx \nu_{p}^{y} (\zeta^{y})
	 P\{\xi^{\Lambda '}_{t}=\zeta^{\Lambda '} \} \\
&	\approx \nu_{p}^{y} (\zeta^{y})\nu_{p}^{\Lambda '} (\zeta^{\Lambda '})\\
&	=\nu_{p}^{\Lambda} (\zeta^{\Lambda }).
\end{split} 
\end{equation*}
\end{proof}

\section{Stationary Measures, Ergodicity, and Mixing}\label{sec:ergodicity}

\subsection{Examples of Stationary Measures}\label{ssec:examples}
We have shown that the product Bernoulli-$p$ measure $\nu_{p}$ is
stationary for the northeast process. This is not, however, the only
stationary measure: for instance, the measure $\nu_{0}$ that assigns
probability one to the configuration with all spins $0$ is also
stationary, and ergodic. In fact, there are infinitely many distinct
(and mutually singular) ergodic stationary distributions. A
denumerable family may be built as follows.

Let $\Gamma$ be an infinite, connected subset of $\zz{Z}^{2}$ such
that 
\begin{enumerate}
\item [(a)] for every site $y\in \Gamma$, at least one of the
southerly or westerly nearest neighbors of $y$ is also an element of
$\Gamma$; and 
\item [(b)] every connected component of $\zz{Z}^{2}\setminus \Gamma$
is finite.
\end{enumerate}
Consider the initial configuration $\xi_{0}$ in which every site in
$\Gamma$ is assigned spin $0$ and every other site spin $1$. If the
northeast process started in this configuration then every site in
$\Gamma$ will remain frozen at spin $0$, because any query tree (see
sec.~\ref{ssec:ppConstruction}) that begins at such a site will have
leaves in $\Gamma$. Similarly, every site $x$ one of whose southerly
or westerly nearest neighbors is in $\Gamma$ will remain frozen at
spin $1$ forever. Now consider the evolution $\xi^{\Lambda}_{t}$ in a
connected component $\Lambda $ of $\zz{Z}^{2}\setminus \Gamma$:
because $\Lambda$ is finite, and because its southwest border consists
of sites that must remain frozen at spin $1$, the process
$\xi^{\Lambda}_{t}$ is an ergodic finite-state Markov process, and so the
distribution of $\xi^{\Lambda}_{t}$ must converge as $t \rightarrow
\infty$ to a stationary distribution $\mu^{\Lambda}$.  The product
$\lambda_{\Gamma}$ of these measures $\mu^{\Lambda}$ over all
components $\Lambda$ with the point mass on the zero configuration in
$\Gamma$ is a stationary distribution for the northeast process. Since
each of the factors $\mu^{\Lambda}$ is stationary and ergodic for the
restricted process $\xi^{\Lambda}_{t}$, the product measure
$\lambda_{\Gamma}$ will be ergodic for the northeast process.  Observe
that the measures $\lambda_{\Gamma}$ and $\lambda_{\Gamma '}$ will be
mutually singular if $\Gamma \not =\Gamma '$, because different sets
of sites ($\Gamma$ and $\Gamma '$) remain frozen forever.

The preceeding construction shows that there are infinitely many
ergodic stationary distributions for the northeast model. None of
these distributions is translation invariant. It is natural to ask
if there exist translation invariant stationary distributions other
than the product Bernoulli measures $\nu_{p}$ and $\nu_{0}$. The
answer is yes, as the following argument shows. Fix an integer $m\geq
2$, and let $\Gamma_{0}$ be the set of sites with at least one
coordinate is divisible by $m$. For each $x= (x_{1},x_{2})\in
\zz{Z}^{2}$ with coordinates satisfying $0\leq x_{i}<m$, let
$\Gamma_{x}=\Gamma_{0}+x$ be the translation of $\Gamma_{0}$ by
$x$. Note that each of the sets $\Gamma_{x}$ is invariant under
translations by $(m,0)$ and $(0,m)$, and that
$\Gamma_{x}+y=\Gamma_{x+y}$ where the addition in the subscript is
done mod $m$. Define
\begin{equation}\label{eq:mixture}
	\mu = \frac{1}{m^{2}}\sum_{x \in \zz{Z}_{m}^{2}}
		\lambda_{\Gamma_{x}}.
\end{equation}
Then $\mu$ is certainly translation invariant, because a translate of
$\Gamma_{x}$ is another $\Gamma_{y}$, and $\mu$ is stationary for the
northeast process, because it is a mixture of stationary
distributions. Note that $m$ is not ergodic, because the initially
frozen set $\Gamma_{x}$ must remain frozen forever.

\subsection{Characterization of Ergodic, Translation Invariant
Measures}\label{ssec:eti} We now show that if $\mu$ is an ergodic,
translation invariant stationary distribution for the northeast
model with flip parameter $p$ then either $\mu =\nu_{p}$ or $\mu
=\nu_{0}$.  Denote by $P_{\mu}=P\times \mu$ the probability measure on
$\Omega \times \{0,1 \}^{\zz{Z}^{2}}$ according to which the initial
configuration $\xi_{0}$ is chosen randomly from $\mu$ and the marked
Poisson processes used in the construction of
sec.~\ref{ssec:ppConstruction} are built on $(\Omega ,\mathcal{F},P)$.
Under $P_{\mu }$, the distribution of $\xi_{0}$ is translation
invariant, and so there exists a constant $0\leq r\leq 1$ such that
$P_{\mu}\{\xi_{0} (x)=1 \}=r$ and hence (by stationarity)
\begin{equation}\label{eq:r}
	P_{\mu}\{\xi_{t} (x)=1 \}=r \quad \forall \, x\in \zz{Z}^{2},
	\, \forall t \geq 0.
\end{equation}
If $r=0$ then $\mu =\nu_{0}$. If $r=1$ then at every rational time,
every site is in spin $1$; but this is impossible unless $p=1$,
because it would imply that no site ever flips to spin $0$, as spin
values are held for time intervals of positive duration. Thus, if
$0<p<1$ then either $r=0$ or $0<r<1$.

Suppose that $0<r<1$. Since $\mu$ is assumed to be ergodic, 
the Birkhoff ergodic theorem implies that the long-time average spin
value at any site $x$ converges to $r$, almost surely. Thus, with
$P_{\mu}-$probability $1$, site $x$ flips its spin infinitely often,
and so $P_{\mu}\{\tau_{x}<\infty \,\forall \, x\}=1$. Hence, Corollary
\ref{corollary:mixing} implies that for $\mu -$almost every initial
configuration $\zeta $ the conditional joint distribution of any
finite spin block $\xi^{\Lambda}_{t}$ given $\xi_{0}=\zeta $ converges
weakly as $t \rightarrow \infty$ to the product Bernoulli measure
$\nu^{\Lambda}_{p}$. It now follows that $r=p$ and $\mu =\nu_{p}$. 

\subsection{Ergodicity and Mixing of $\nu_{p}$ for
$p>1-\beta_{c}$}\label{ssec:aboveCritical} To complete the proof of
Theorem~\ref{theorem:1} and the first assertion of Theorem
\ref{theorem:2} we must show that if $p\geq 1-\beta_{c}$ then
$\nu_{p}$ is ergodic and mixing for the northeast model. The key is that if
$p\geq 1-\beta_{c}$ then for $\nu_{p}-$almost every spin configuration
$\zeta$ there are no infinite southwest clusters of $0$s. In
particular, the size of the cluster containing the origin is almost
surely finite.

\begin{lemma}\label{lemma:onOff}
Define $M^{x}_{t}$ to be the number of flip opportunities at site $x$
up to time $t$, and denote by $x_{S}, x_{W}$ the southerly and
westerly nearest neighbors of $x$. Then  for every initial
configuration $\zeta$,
\begin{equation}\label{eq:onOff}
	P (M^{x} (\infty)<\infty \quad \text{and}\quad 
	M^{x_{S}} (\infty)=M^{x_{W}} (\infty)=\infty \, |
	\,\xi_{0}=\zeta)
	=0.
\end{equation}
\end{lemma}

\begin{proof}
It is enough to show that almost surely on the event
\[
G:=\{M^{x_{S}}
(\infty)=M^{x_{W}} (\infty)=\infty \}
\]
the spins at sites $x_{S}$ and
$x_{W}$ will both be $1$ at arbitrarily large times $t$. Denote by
$\tau^{k}_{y}$ the time of the $k$th flip opportunity at site $y$. On
the event $G$ it must be that $\tau^{k}_{y}<\infty$ for all $k\geq 1$
and both $y=x_{S},x_{W}$.  At each time $\tau^{k}_{y}$, the spin at
site $y$ is reset by a toss of a $p-$coin. Since the spin values at
$y=x_{S}$ and $y=x_{W}$ play no role in determining the reset times at
either $x_{S}$ or $x_{W}$, it follows that the \emph{pair} of spins at
$x_{S},x_{W}$ will be reset by independent $p-$coin tosses infinitely
many times. Hence, the spins will both be $1$  at indefinitely large
times, w.p.1.
\end{proof}

Assume now that the initial configuration $\xi_{0}$ is chosen at
random according to the product Bernoulli-$p$ measure $\nu_{p}$,
independently of the marked Poisson processes used to determine the
time evolution, and assume that $p\geq 1-\beta_{c}$. To show that the
process $\xi_{t}$ is mixing (and therefore also ergodic), it suffices,
by Corollary \ref{corollary:mixing}, to show that $\tau_{x}<\infty$
almost surely for each site $x$. We will show that in fact $M^{x}
(\infty)=\infty$ a.s.

Suppose to the contrary that for some site $x$ there is positive
probability, say $\varrho >0$, that $M^{x}(\infty)<\infty $. Then by
Lemma \ref{lemma:onOff}, for at least one of the sites $y=x_{S},x_{W}$
it must be the case that $M^{y} (\infty)<\infty$. It then follows that
for at least one of these two sites, the spin at this site stabilizes
at $0$, because if $\xi_{t} (y)=1$ eventually and the other southwesterly
nearest neighbor $y'$ has infinitely many flip opportunities, then
$\xi_{t} (y')=1$ for arbitrarily large times $t$ and so $x$ must have
infinitely many flip opportunities, contradicting the hypothesis that
$M^{x} (\infty)<\infty$. By induction, it follows that there is an
infinite sequence $y_{n}$ of sites such that (i) each $y_{n+1}$ is
either the southerly or westerly nearest neighbor of $y_{n}$; and (ii)
the spin at each $y_{n}$ eventually stabilizes at $0$.

This contradicts the hypothesis that $p\geq 1-\beta_{c}$: In
particular, it implies (by translation invariance) that for each
$K<\infty $ there is probability at least $\varrho $ that at large
times $t$ the origin will belong to a size-$\geq K$ southwest cluster
of sites all with spin $0$. Since by stationarity of $\nu_{p}$ the
size of the southwest cluster of $0$s containing the origin has the
same distribution at time $t$ as at time $0$, this is impossible.
\qed

\section{Exponential Decay of Correlations}\label{sec:expMixing}

In this section we shall establish the exponential decay of
correlations (\ref{eq:expMixing}) for all parameter values $p>1/2$. 
Without loss of generality, we may assume that the finite set
$\Lambda$  in (\ref{eq:expMixing}) is a square, since every finite set
is contained in a square; and furthermore, we may assume that the
northwest corner of the square $\Lambda$ is the origin
$(0,0)$. 

Denote by $\partial_{-}\Lambda$ the exterior southwest boundary of
$\Lambda$, that is, the set of all sites not in $\Lambda$ that border
$\Lambda $ to the south or west. At any (Markov) time when the spins
at all sites in $\partial_{-}\Lambda$ are $1$ it is possible for the
sites in $\Lambda $ to begin to reset, starting from the southwest
corner and proceeding north and east. In particular, it is possible
that all spins in $\Lambda $ will flip to $1$ \emph{before} any of the
spins in $\partial_{-}\Lambda$ reset. Moreover, at any (Markov) time
when all of the spins in $\Lambda$ are $1$, it is possible that the
spins will then reset one at a time  starting from the
northeast corner and proceeding to the southwest, in the order
indicated below (for a $3\times 3$ square):
\[
\begin{matrix}{}
4 & 2 & 1 \\
7 & 5 & 3 \\
9 & 8 & 6
\end{matrix}
\]
Call such an occurrence a
(total) $\Lambda -$\emph{reset}. Define
\begin{equation}\label{eq:tauLambda}
	\tau_{\Lambda}= \text{time of first} \;\Lambda\text{--reset}.
\end{equation}
Observe that at time $\tau_{\Lambda}$ the spin at the southwest corner
of $\Lambda $ has just reset. At this time all spins in $\Lambda$ have
been reset at least once, in order starting from the northeast corner,
and so the conditional distribution of the configuration
$\xi^{\Lambda}_{\tau_{\Lambda}}$, given the evolution to time
$\tau_{\Lambda}$ of the process in the region southwest of $\Lambda$,
is the Bernoulli product measure $\nu^{\Lambda}_{p}$.

\begin{lemma}\label{lemma:tauLambda}
To prove (\ref{eq:expMixing}) it suffices to prove that  $\forall t\geq 0$,
\begin{equation}\label{eq:tL}
	 P^{\nu_{p}}_{p}\{\tau_{\Lambda}\geq t \}
	\leq C_{\Lambda }\exp \{-\alpha_{\Lambda}t \}.
\end{equation}
\end{lemma}

\begin{proof}
This follows from Proposition \ref{proposition:resets} by a similar
argument as in the proof of Corollary \ref{corollary:mixing}.
\end{proof}

Whether or not a $\Lambda -$reset will occur quickly depends on the
current values of the spins at the sites to the southwest of
$\Lambda$, in particular, on those that directly border on
$\Lambda$. Since any such site with spin value $0$ will freeze the
spins at its northern and eastern neighbors, it is impossible for
$\Lambda $ to reset until all the southwest boundary sites have
assumed spin value $1$, unless all spins in $\Lambda$ are initially at
spin $1$. Similarly, any of the boundary sites initially at spin $0$
cannot flip until their southwest neighbors have assumed spin value
$1$, and so on. Thus, at any time $t$ the obstruction to a $\Lambda
-$reset is the union of the oriented $0-$clusters attached to the
southwest boundary of $\Lambda$. Denote this union by
$K^{\Lambda}_{t}$: that is, $K^{\Lambda}_{t}$ is the set of all sites
$y$ for which there exists a (northwest oriented) path of $0$s in the
configuration $\xi_{t}$ beginning at $y$ and ending at the southwest
boundary of $\Lambda$.  Fix $t_{*}>0$, and define stopping times
\[
	T_{0} < S_{0} < T_{1} < S_{1} < \dotsb 
\]
by setting
\begin{align*}
T_{0}&= 0;\\
S_{n}&= \min \{t\geq T_{n}\, : \, K^{\Lambda}_{t}=\emptyset \};\\
T_{n}&= \min \{t>S_{n-1}\, : \, K^{\Lambda}_{t}\not =\emptyset \}\wedge t_{*},
\end{align*}
where $t_{*}$ is a fixed, nonrandom time. 

\begin{lemma}\label{lemma:clusterProcess}
Assume that $p>2/3$. Then there exist constants $C,\alpha >0$,
possibly depending on the size of the square $\Lambda$ and on the
choice of the  constant $t_{*}$, such
that for all $n$ and all $t>0$,
\begin{equation}\label{eq:clusterProcess}
	P (S_{n}-T_{n}\geq t \, | \, \mathcal{I}^{\infty}_{T_{n}})
	\leq C \exp \{-\alpha t \}.
\end{equation}
\end{lemma}

\begin{proof}
Let $X_{t}:=|K^{\Lambda}_{t}|$ denote the cardinality of the
$0-$cluster set attached to $\Lambda$ at time $t$; then $S_{n}$ is the
first time after $T_{n}$ when $X_{t}=0$.  At any time when the cluster
$K^{\Lambda}_{t}$ is nonempty, its exterior southwest boundary
consists of $B\geq 1$ sites, all at spin $1$, and its interior
southwest boundary (the set of sites that are
flip-eligible) consists of $A\geq 1$ sites, all at spin $0$; since
each interior boundary site is bordered by at most two exterior
boundary sites, 
\begin{equation}\label{eq:ABRatio}
	B \leq 2A.
\end{equation}
The cluster can grow or shrink by only one site at a time, either by
addition of a site in the exterior southwest boundary or deletion of a
site in the interior southwest boundary. Since there is always at least
one interior boundary site that is flip-eligible, the rate
at which jumps in $X_{t}$ occur is at least $p$. Moreover,
at any addition/deletion event, the conditional probability that the
event is a deletion, given the current configuration, is at least
\[
	pA/ (pA+ (1-p)B):
\]
this is because the instantaneous rate of deletions is $pA$, whereas
the instantaneous rate of additions is \emph{at most} $(1-p)B$. (All
of the spins in the exterior southwest boundary of $K^{\Lambda}_{t}$
are $1$, but some of these may be frozen.) Since $B\leq 2A$ and
$p>2/3$, 
\[
	pA/ (pA + (1-p)B)>1/2.
\]
Thus, the jump process $X_{t}$ is stochastically dominated by a
reflecting random walk on the nonnegative integers with negative
drift. It follows routinely that the first  passage time to $0$ has an
exponentially decaying tail.
\end{proof}

\begin{proof}
[Proof of (\ref{eq:tL})] Let $G_{n}$ be the event that there is a
$\Lambda -$reset between times $S_{n}$ and $T_{n+1}$. Since at time
$S_{n}$ the $0-$cluster set $K^{\Lambda}_{S_{n}}$ is empty, the
conditional probability of $G_{n}$, given the history of the process
$\xi_{t}$ up to time $S_{n}$, is bounded below by a constant $\delta
>0$. (Recall that the spins in $\partial_{-}\Lambda$ must stay fixed
at $1$ until there is an occurrence in one of the Poisson processes
$N^{y}_{t}$ attached to a site $y\in \partial_{-}\Lambda$, regardless
of what happens elsewhere.) Hence, by Lemma
\ref{lemma:clusterProcess}, 
\begin{equation}\label{eq:iidBound}
	\tau_{\Lambda}\leq S_{0}+\sum_{j=1}^{H}R_{j}
\end{equation}
where $S_{0}, H$, and $R_{1},R_{2},\dotsc$ are mutually independent,
$H$ is the number of tosses of a $\delta -$coin until the first Head,
and
\[
	P\{R_{j}\geq t \}=C \exp \{-\alpha t \} \quad 
	\forall \; t\geq \log C/\alpha .
\]
Elementary arguments show that the random variable on the right side
of (\ref{eq:iidBound}) has an exponentially decaying tail.
\end{proof}

\section{Influence Propagation}\label{sec:shape}

A key component of the Aldous-Diaconis argument (although not
explicitly noted as such) is a proof that the influence of a spin $1$
is propagated forward at a definite linear rate. In two dimensions,
propagation of influence is a much more subtle business, as there is
no longer a unique path along which it occurs. A natural way to
measure influence propagation in the northeast model is by use of
\emph{influence regions} $R_{t}$ defined as
follows. Consider an initial configuration  in
which all sites in the first quadrant are initially set to spin
$0$, and all other sites are initially set to random (by independent
$p-$coin tosses).
For each time $t\geq 0$ define 
\begin{equation}\label{eq:influencedRegion}
	R_{t}=\{\text{sites in quadrant 1 that are flipped at least
once by time} \;t \}.
\end{equation}

Figure $1$ below shows the results of a simulation (done by
\textsc{Marc Coram} of the University of Chicago) in which sites in
quadrants 2, 3, and 4 are initially determined by independent $p-$coin
tosses, and $p=.8$. The figure shows the state of the system at time
$t=1044.65$; the size of the box shown is $500\times 500$. The
influence region $R_{t}$ is shown in white, and the blue region
$Q_{t}$ consists of those sites that have been ``queried'' at least
once by time $t$. See section \ref{sec:bg} for an explanation of the
query process. The black specks indicate those sites in $R_{t}$ or
$Q_{t}$ whose spin values are $0$ at time $t$.

\begin{figure} 
\label{NorthEast}
 \includegraphics[width=3.0in,height=3.0in]{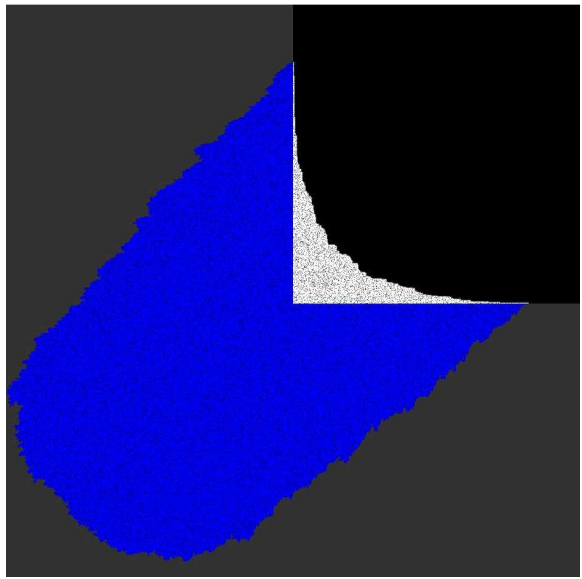}
 \caption{Influence Propagation in Northeast Model} 
  \end{figure}

The results of this and
other simulations suggest the following conjecture.

\begin{conjecture}\label{conjecture:1}
For $p>p_{c}$,
the influence region $R_{t}$ as a definite, nonrandom \emph{limit
shape} $S$. That is, with $\hat{R}$ denoting the union  of all unit
squares centered at points of $R$,
\begin{equation}\label{eq:shapeTheorem}
	\frac{\hat{R}_{t}}{t}
	\longrightarrow S
\end{equation}
in the Hausdorff metric, for some nonrandom set $S$ contained in
the first quadrant of $\zz{R}^{2}$.
\end{conjecture}

Note that because the Northeast proces is not additive (in the  sense
of \cite{liggett}), the usual methods for proving shape theorems (see,
e.g., \cite{cox} and \cite{richard}) are
not applicable.

%%%%%%%%%%%%%%%%%%%%%%%%%%%%%%%%%%%%%%%%%%%%%%%%%%%%%%%%%%%%%%%%%%%%%%

\end{document}